# A NONPARAMETRIC APPROACH TO THE ESTIMATION OF LENGTHS AND SURFACE AREAS

By Antonio Cuevas,[1] Ricardo Fraiman[1] and
Alberto Rodríguez-Casal[2]

*Universidad Autónoma de Madrid, Universidad de San Andrés
and Universidad de Santiago de Compostela*

The Minkowski content $L_0(G)$ of a body $G \subset \mathbb{R}^d$ represents the boundary length (for $d=2$) or the surface area (for $d=3$) of $G$. A method for estimating $L_0(G)$ is proposed. It relies on a nonparametric estimator based on the information provided by a random sample (taken on a rectangle containing $G$) in which we are able to identify whether every point is inside or outside $G$. Some theoretical properties concerning strong consistency, $L_1$-error and convergence rates are obtained. A practical application to a problem of image analysis in cardiology is discussed in some detail. A brief simulation study is provided.

**1. Introduction.** The estimation of the surface area of a body $G$ in the Euclidean space $\mathbb{R}^d$ ("surface area" amounts to "boundary length" in the bidimensional case $d=2$) has been extensively considered in stereology; see [1, 2, 12]. We are concerned here with this problem from a different point of view, using the approach and tools of nonparametric statistics and, more specifically, of nonparametric set estimation; see, for example, [6] for a survey.

In a way, the length and surface area estimation problem can be seen as a further, more difficult, stage in set estimation theory, after the early developments concerned with the estimation of volume (associated with the $L_1$ (measure) distance; see [8]), "visual" shape (associated with the Hausdorff metric; see [5]), level sets [3, 13, 16, 19, 20] and boundaries [7]. We will

Received June 2005; revised March 2006.
[1]Supported in part by Spanish Grant MTM2004-00098.
[2]Supported in part by Spanish Grants MTM2005-00820 and PGIDIT06PXIB207009PR.

see, in fact, that, while the sample data in nonparametric set estimation theory comes usually from random points selected inside the set of interest, $G$, we will need here additional information given by sample points coming from outside $G$ (see the beginning of Section 2). The estimation of boundary length has also some practical interest. For example, in medical imaging the boundary length appears in connection with the notion of "Contour Index" (see, e.g., [11]), a shape measurement used as an auxiliary diagnostic criterion. These ideas are developed in more detail in Section 4.

At this point it might be useful to point out what we mean by "nonparametric approach" in order to clarify its main differences with the stereological point of view for these problems:

(a) Unlike the stereological approach, we are not concerned with unbiased estimation, but with asymptotic properties such as consistency and convergence rates.

(b) The proposed estimator is intended to work asymptotically in any dimension $d$ under quite general shape restrictions. It depends on a smoothing parameter which must be carefully chosen.

(c) Our method will provide as a by-product an estimator of the boundary of the body $G$ under study. In contrast, stereological methods are not usually concerned with the global estimation of sets; they are rather focused on the estimation of some real parameter (length, volume, surface area, ...).

(d) The sample data consists of randomly selected points. In stereology the available information for estimating lengths and surface areas usually comes either from one- or two-dimensional sections or from systematic grids.

Our aim is to obtain an easy-to-implement automatic method valid for the analysis of a wide class of images. As a first step we should clearly establish what we mean by "surface area." The Hausdorff measure (see, e.g., [14]) provides a suitable general definition of this concept. This definition, however, is not always very convenient from the point of view of mathematical handling and effective evaluation. So we will use instead the following simpler, less general notion (which coincides with the Hausdorff measure, up to a constant factor, in regular cases): The surface area of a body $G \subset \mathbb{R}^d$ is given by the *Minkowski content* (see [14], Chapter 2),

$$(1) \qquad L_0(G) = \lim_{\varepsilon \to 0} \frac{\mu(B(\partial G, \varepsilon))}{2\varepsilon},$$

provided that this limit exists and it is finite. Here $\mu$ stands for the ordinary Lebesgue measure on $\mathbb{R}^d$, $\partial G$ denotes the boundary of $G$ and, for any $A \subset \mathbb{R}^d$, $B(A, \varepsilon)$ is the "outer parallel set" $B(A, \varepsilon) := \bigcup_{x \in A} B(x, \varepsilon)$, where $B(x, \varepsilon)$ denotes the closed ball with center $x$ and radius $\varepsilon$. While the Minkowski content fails to satisfy some interesting properties, such as $\sigma$-additivity, it has a clear intuitive basis and is sufficient for most practical purposes.

This paper is organized as follows. The estimator is introduced in Section 2. Its basic statistical properties concerning asymptotic behavior, bias and variability are established in Section 3. A real-data application in cardiology is discussed in Section 4. A brief Monte Carlo study is presented in Section 5. Section 6 is devoted to the proofs.

**2. The sampling model and the proposed estimator.** Let $G \subset \mathbb{R}^d$ be a body whose Minkowski content $L_0 = L_0(G)$ is well defined, strictly positive and finite. Our goal is estimating $L_0$ which for $d = 2$ represents the boundary length and for $d = 3$ the surface area. Without loss of generality, we will assume that $G$ is a subset of the open unit square $(0,1)^d$.

The sampling information is given by i.i.d. observations $(Z_1, \delta_1), \ldots, (Z_n, \delta_n)$ of a random variable $(Z, \delta)$, where $Z$ is uniformly distributed on the unit square $[0,1]^d$ and $\delta = 1$ if $Z \in G$, $\delta = 0$ if $Z \notin G$. This means that, with probability one, given a sample of points on the unit square, we are able to decide whether or not they belong to the "green area" $G$ or to the "red" one, $R = [0,1]^d \setminus G$.

It will be convenient to use the following notation. Let us denote by $P_X$ and $P_Y$ the conditional distributions of the "green" and "red" observations, that is, the distributions of $Z|\{\delta = 1\}$ and $Z|\{\delta = 0\}$. Observe that $P_X$ and $P_Y$ are both uniform on $G$ and $R$, respectively. Now, given $z \in [0,1]^d$ and $\varepsilon \geq 0$, denote by $G_z(\varepsilon)$ and $R_z(\varepsilon)$, respectively, the numbers of green and red sample observations belonging to the ball $B(z, \varepsilon)$, that is,

$$
\begin{aligned}
G_z(\varepsilon) &\equiv G_{n,z}(\varepsilon) = \sum_{i=1}^{n} \mathbb{I}_{\{\delta_i = 1, \|Z_i - z\| \leq \varepsilon\}}, \\
R_z(\varepsilon) &\equiv R_{n,z}(\varepsilon) = \sum_{i=1}^{n} \mathbb{I}_{\{\delta_i = 0, \|Z_i - z\| \leq \varepsilon\}}.
\end{aligned}
\tag{2}
$$

Clearly, $G_z(\varepsilon)$ has a binomial distribution with parameters $n$ and $p_X(z, \varepsilon) = P(\|Z - z\| \leq \varepsilon, \delta = 1) = \mu(G) P_X(B(z, \varepsilon))$. Similarly, $R_z(\varepsilon)$ has a binomial distribution with parameters $n$ and $p_Y(z, \varepsilon) = (1 - \mu(G)) P_Y(B(z, \varepsilon))$.

Let $\{\varepsilon_n\}$ be a deterministic sequence of positive numbers which converges to zero as $n$ tends to infinity. Denote $T = \partial G$. We propose the following estimator for the "dilated boundary," $B(T, \varepsilon_n)$:

$$
T_n = \{z \in [0,1]^d : R_z(\varepsilon_n) \geq 1 \text{ and } G_z(\varepsilon_n) \geq 1\}.
\tag{3}
$$

The simple intuitive idea behind $T_n$ is to consider those points $z$ in whose vicinity green and red points coexist. Of course, we could "robustify" this estimator by replacing the condition $R_z(\varepsilon_n) \geq 1$ and $G_z(\varepsilon_n) \geq 1$ with $R_z(\varepsilon_n) \geq r_1$ and $G_z(\varepsilon_n) \geq g_1$, for some fixed integer numbers $r_1 > 1$ and $g_1 > 1$. This modified estimator (which will not be considered here) would be smoother and less noisy than the original version (3) at the expense of some efficiency loss.

Finally, the definition (1) for $T_n$ suggests the following natural estimator for $L_0 = L_0(G)$:

$$L_n = \frac{\mu(T_n)}{2\varepsilon_n}. \tag{4}$$

As usual, the nonparametric estimator (4) depends on a smoothing parameter $\varepsilon_n$ which must be carefully chosen. In general, it should tend to zero slowly enough. The theoretical results of Section 3 will provide some additional insights in this respect.

Note that the proposed method could be useful even in those cases where the image $G$ is completely known (e.g., we could have a picture of $G$), but it is too complicated for directly measuring its boundary. Then the sample $Z_1, \ldots, Z_n$ can be artificially generated provided that we are able to decide whether $Z_i$ belongs to $G$ or not. So, in some sense, (4) can be seen as a "stochastic" algorithm to approximate $L_0$. This idea will be further developed in Section 4.

**3. Theoretical results.** We analyze in this section the properties of the estimator $L_n$ of the Minkowski content, $L_0 = L_0(G)$.

3.1. *Strong consistency.* The almost sure (a.s.) convergence of $L_n$ to $L_0$ is established in Theorem 1 below. The "standardness" hypothesis (a) prevents the set $G$ from having "too sharp" inlets and peaks along the boundary $T$. This condition has been previously used in set estimation (see, e.g., [7]).

THEOREM 1. *Let us assume the following conditions.*

(a) *The sets $G$ and $R$ are both standard in $T$, that is, there exists a constant $C > 0$ such that, for small enough $\varepsilon$,*

$$P_X(B(t,\varepsilon)) \geq C\mu(B(t,\varepsilon)) \quad and \quad P_Y(B(t,\varepsilon)) \geq C\mu(B(t,\varepsilon)) \qquad for\ all\ t \in T.$$

(b) *The sequence $\{\varepsilon_n\}$ satisfies*

$$\varepsilon_n \to 0 \quad and \quad \frac{n\varepsilon_n^d}{\log n} \to \infty.$$

*Then*

$$L_n = \frac{\mu(T_n)}{2\varepsilon_n} \to L_0, \qquad a.s.$$

Observe that the conditions imposed in (b) on the sequence $\varepsilon_n$ of smoothing parameters are identical to those required for the strong consistency of kernel density estimators (see, e.g., [15]). However, as we will see below, the role of the smoothing parameter is quite different in both setups.

3.2. *The function $L(\varepsilon)$.* For a given value of $n$, the estimator $L_n$ provides, in fact, an estimation for $L(\varepsilon_n) := \mu(B(T, \varepsilon_n))/(2\varepsilon_n)$ which, in turn, is an approximation of the target value $L_0$. Thus, in order to assess the accuracy of the estimator $L_n$, it is interesting to get more precise information on the difference $|L(\varepsilon) - L_0|$. We next show that, under some smoothness assumptions, $L(\varepsilon)$ is differentiable at $\varepsilon = 0$, which entails $|L(\varepsilon) - L_0| = O(\varepsilon)$. Indeed, note that

$$B(T, \varepsilon) = B(G, \varepsilon) \cap B(R, \varepsilon),$$

which leads to

(5) $\qquad \mu(B(T,\varepsilon)) = \mu(B(G,\varepsilon)) + \mu(B(R,\varepsilon)) - \mu(B([0,1]^d, \varepsilon)).$

Thus, the point is to have some idea about the structure of the "dilated measures" on the right-hand side of (5), when considered as functions of $\varepsilon$. If $G$ is assumed to be convex, the classical Steiner formula (see, e.g., [17], page 197) establishes that $\mu(B(G, \varepsilon))$ is a polynomial in $\varepsilon$ of degree at most $d$. Unfortunately, this result is not useful in our case, as the hypothesis of convexity for $G$ could be too restrictive (e.g., in image analysis) and, in any case, it cannot be assumed simultaneously for both $G$ and $R = [0,1]^d \setminus G$, except in trivial situations. However, we will be able to prove the required differentiability property for $L(\varepsilon)$ by combining some ideas of mathematical morphology (which we will use to impose the appropriate regularity conditions on $G$) with a (partial) generalization of Steiner's formula proved by Federer [10]. He imposes a *positive reach condition* closely related to the following *rolling condition* often used in set estimation (see, e.g., Walther [20]): It is said that a ball *can roll along $T = \partial G$ outside $G \subset \mathbb{R}^d$* if there exists $r_0 > 0$ such that, for all $r \leq r_0$ and $x \in T$, there exists a closed ball of radius $r$, $B_x$, such that $B_x \cap G = \{x\}$.

A deep study of this *outer rolling condition*, including some interesting equivalences, is due to Walther [21], Theorem 1. This condition arises in *mathematical morphology*, a branch of the huge current theory of image analysis; see [18]. It has also appeared, under a slightly different form, in contexts not directly related to image analysis. In a similar vein, Federer [10] defines the *reach of $G$* as the largest (possibly $\infty$) value $r_0$ such that if $x \in \mathbb{R}^d$ and the distance from $x$ to $G$ is smaller than $r_0$, then $G$ contains a unique point nearest to $x$. For our purposes of better understanding the nature of the function $L(\varepsilon)$, it will be particularly useful to employ a generalization of Steiner's formula obtained by Federer ([10], Theorem 5.6). This result establishes that, for any set $G$ of positive reach $r_0$, the function $\mu(G, \varepsilon)$ coincides locally [for $\varepsilon \in (0, r_0)$] with a polynomial of degree at most $d$ whose independent term is $\mu(G)$.

Thus, if we assume that both $G$ and $R$ satisfy the positive reach condition we may use Federer's theorem, together with (5), to conclude that $\mu(T, \varepsilon)$ coincides in the interval $(0, r_0)$ with $P(\varepsilon)$, where $P$ denotes a polynomial of degree at most $d$ with a null independent term. Note that (by the assumption made on the finiteness of the Minkowski content $L_0$) the coefficient of $\varepsilon$ in $P(\varepsilon)$ must necessarily coincide with $2L_0$ so that $L(\varepsilon) - L_0$ is a polynomial in $\varepsilon$ with a null independent term. In particular, we get that $L(\varepsilon)$ is differentiable at $\varepsilon = 0$.

3.3. *Bounds for $E(L_n)$: $L_1$-consistency and convergence rates, variability and bias.* It is not hard to show (see the proof of Statement 1 in the proof of Theorem 1) that, with probability one, $T_n \subset B(T, \varepsilon_n)$ and, therefore,

(6) $$L_n \leq L(\varepsilon_n) \quad \text{a.s.}$$

This means that $L_n$ tends to underestimate $L_0$ for those "regular" sets where the values of the function $L(\varepsilon) = \mu(B(T,\varepsilon))/(2\varepsilon)$ are very close to $L(0) := L_0$ for small values of $\varepsilon$. In the bidimensional case the simplest example is given by the circle, for which $L(\varepsilon) \equiv L_0$.

The following result provides a lower bound for $E(L_n)$.

THEOREM 2. *Assume that the standardness condition* (a) *in Theorem 1 holds. Assume also that the function $F(\varepsilon) := \mu(B(T,\varepsilon))$ is differentiable in a neighborhood of $0$ and the derivative $F'$ is continuous at $0$. Then*

(7) $$E(L_n) \geq L(\varepsilon_n) - I_n,$$

*where $I_n = \frac{1}{\varepsilon_n} \int_{B(T,\varepsilon_n)} \exp(-Kn(\varepsilon_n - d(z,T))^d) \, dz$, $K$ being a positive constant and $d(z,T) = \inf\{\|z - t\| : t \in T\}$. Also,*

(8) $$I_n = O((n\varepsilon_n^d)^{-1/d}).$$

The proof is given in Section 6. Note that, according to the discussion in Section 3.2, if we assume that both $G$ and $R$ fulfill the positive reach property, then the function $F(\varepsilon) = \mu(B(T,\varepsilon))$ coincides in a neighborhood of $0$ with a polynomial of degree $\leq d$, so it is certainly differentiable at $0$ with a continuous derivative.

The following corollary (the proof is in Section 6) provides a condition for the $L_1$-consistency, as well as an upper bound for the $L_1$-convergence rate of the estimator $L_n$.

COROLLARY 1. (a) *Under the same conditions of Theorem 2, we have*

(9) $$E|L_n - L_0| \leq I_n + |L(\varepsilon_n) - L_0|.$$

*As a consequence, the standard conditions for consistency, $\varepsilon_n \to 0$ and $n\varepsilon_n^d \to \infty$, are also sufficient here to ensure the $L_1$-consistency $E|L_n - L_0| \to 0$.*

(b) *By assuming further that $G$ and $R$ satisfy the positive reach condition mentioned in Section 3.2, we have that the optimal order for the bound (9) is $O(n^{-1/2d})$, which is attained for $\varepsilon_n = n^{-1/2d}$.*

Not surprisingly, the bound $O(n^{-1/2d})$ corresponds to a rather slow convergence rate. We do not believe that the exact rate can improve much on this bound. Recall that the typical rates for the much easier problem of consistently estimating the boundary $\partial G$, with respect to the Hausdorff metric, are of type $O((\log n/n)^{1/d})$ [7] even under the assumption of convexity on $G$ [9]. Anyway, in some applications (see Section 4) the estimator $L_n$ is based on artificial (Monte Carlo) samples and the slow convergence rate is not so crucial a problem, as the sample size can, in principle, be increased as much as necessary.

As a further consequence of (6)–(8) we get [under the regularity assumptions imposed in Corollary 1(b)] the following bounds for the $L_1$-variability and the bias:

$$
\begin{aligned}
E|L_n - E(L_n)| &\leq E|L_n - L(\varepsilon_n)| + |L(\varepsilon_n) - E(L_n)| \\
&\leq 2I_n = O((n\varepsilon_n^d)^{-1/d}),
\end{aligned}
\tag{10}
$$

$$
\begin{aligned}
L_0 - E(L_n) &= (L(\varepsilon_n) - E(L_n)) + (L_0 - L(\varepsilon_n)) \\
&= O((n\varepsilon_n^d)^{-1/d}) + O(\varepsilon_n).
\end{aligned}
\tag{11}
$$

Thus, the assumption $n\varepsilon_n^d \to \infty$ guarantees the convergence to zero of the variability around the mean, $E|L_n - E(L_n)|$. Note that this condition is identical to that imposed in the classical ($L_2$ or $L_1$) theory of density estimation in order to control the variability term. However, expression (11) shows that $n\varepsilon_n^d \to \infty$ is also useful to make the bias term tend to zero. This is in sharp contrast with the typical situation in nonparametric functional estimation where $\varepsilon_n \to 0$ usually suffices to kill the bias. The situation here is a bit different: We do need the condition $\varepsilon_n \to 0$, but if the convergence is too fast, the estimator $L_n$ will be biased, underestimating the value of $L_0$. Thus, the bias is also controlled by the condition $n\varepsilon_n^d \to \infty$ which is used in the proof of Theorem 1 to prevent the boundary estimator $T_n$ from having spurious "holes" [that would lead to underestimation of $\mu(B(T, \varepsilon_n))$ by $\mu(T_n)$].

Let us also note that it is interesting to assess the magnitude of the "effective bias" $E(L_n) - L(\varepsilon_n)$. This is particularly useful in practical applications (see Section 4 below) when one is willing to consider $L(\varepsilon_n)$ as a reasonable approximation for $L_0$, thus, accepting a systematic bias which hopefully would affect in a similar way all the images under study. In these cases the focus is on the differences $E(L_n) - L(\varepsilon_n)$, analyzed above.

**4. Applications to image analysis.** Let us first emphasize that our approach is basically aimed at those cases where only partial (random) information is available, rather than dealing with completely known images.

These usually come in a digitized form, but the digitization process is itself an approximation involving nontrivial problems, largely beyond the scope of this paper. The classical book by Serra [18], pages 211–224, provides some deep insights in this regard. Anyway, if we have a "known" image, either in a digitized version or in a "exact" format (e.g., the area inside a known closed curve: see Section 5), it is tempting to check the behavior of the estimator (4), based on Monte Carlo random samples, when used to approximate the boundary length.

In this section we develop this idea and apply it to a medical example.

4.1. *The contour index. A case study in cardiology.* The irregularity in the border of a tumor or an infarcted area is an auxiliary diagnostic criterion for malignancy assessment. The so-called "contour index" (CI) provides a size-independent quantitative measurement of boundary roughness. It is defined (for the case of bidimensional images) as the quotient *boundary length*$/\sqrt{area}$. Its minimal value ($2\sqrt{\pi}$) is attained by the circle. The CI has been used in oncology (see, e.g., [4]) and cardiology [11], but the interpretation of this index in the two scientific fields is somewhat different. A high value of the CI in a tumor usually suggests a high dissemination capacity of the injured area. On the contrary, in cardiology the prognosis of an infarction tends to be worse when the damage is highly concentrated with a "regular" border (which will provide a small CI) rather than disseminated in many small irregular patches.

In order to assess the applicability of our estimation method to real examples, we have analyzed an image (Figure 1, left) of the infarcted heart of a pig. It corresponds to one side of a transversal section of the heart which has been exposed to a histochemical reaction that dyes the living cells. Thus, the infarcted cells fail to catch the color, appearing as a white-grayish area in the upper-right side of the image. This area should not be confounded with the endocardial endothelium (which covers the inner part of the heart). It appears in deep white at the centre of the image. In fact, most of this endothelium white area is not placed in the same plane as the considered transversal section. The jpg file of the original color image (Figure 1, top left) has been digitized in an array of $495 \times 710$ pixels. The information stored in every pixel consists of a vector $(x_1, x_2, x_3)$ indicating the level (on the scale 0–255) of primary colors (red, green and blue) at that point. So, if we consider the position coordinates, every pixel is, in fact, a five-dimensional observation.

4.2. *A stochastic algorithm for calculating the CI.* In the example considered the goal is to identify the infarcted area and give an approximate value for the CI. Our estimation method has been used with the following steps:

1. *Image identification and cleaning.* The image of interest (Figure 1, top left) must be treated in order to clearly decide the precise shape of the infarcted area (a bit blurred in the original picture) whose boundary is to be measured. The problem is to decide whether or not a pixel in the picture corresponds to the infarcted area. We have done this using the classical Fisher linear discriminant function. To put this in more precise terms, two large samples of pixels have been taken in the infarcted and in the noninfarcted area. Then the classical linear discrimination method was applied to classify the remaining points. The classification error was negligible except for the points in the white endothelial area at the centre of the original image (that tended to be confounded with the infarcted cells), where the error rate was appreciable. The result of this automatic discrimination-based treatment is shown in Figure 1 (top right) where the infarcted area has been colored in black but there are also some patches of obviously misclassified endothelial tissue. Thus, a final "manual cleaning" was made to remove these patches. The result is given in Figure 1 (bottom). This was the final image (600×600 pixels) used for the quantitative analysis described in Section 4.3. Let us note that the classification algorithm has been based only on the "color" coordinates of every point. We have disregarded the information provided by the point positions because the use of a linear discrimination method looked particularly unsuitable when these variables are involved.

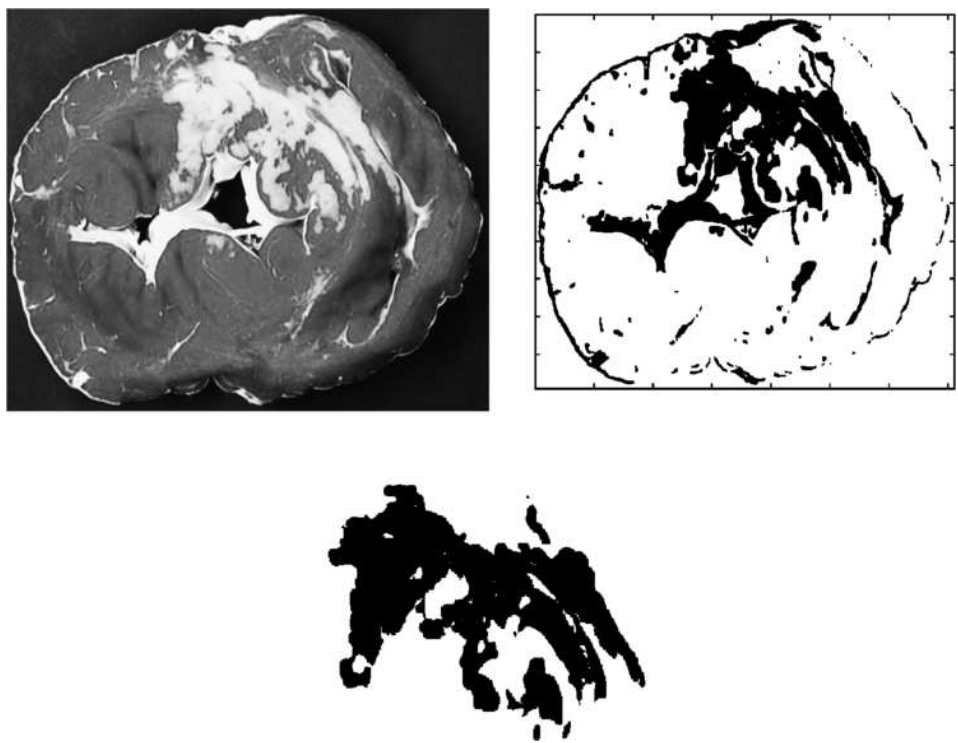

FIG. 1. *An infarcted heart* (top left). *The estimated infarct area* (top right). *The "cleaned" infarct area* (bottom).

By the way, this application of discriminant methods in image cleaning shows the interest in studying discriminant theory from the point of view of image analysis; this would amount to incorporating classification criteria based on shape preservation (connectedness, smoothness), in addition to the usual notions relying on misclassification probabilities.

2. *Monte Carlo sampling and classification.* A large artificial uniform sample $Z_1, \ldots, Z_n$ is drawn in $[0,1]^d$. The classification variable $\delta_i$ is obtained for each $Z_i$: $\delta_i = 1$ when $Z_i$ belongs to the infarcted area, $\delta_i = 0$ otherwise.

3. *Estimation.* As indicated in Section 3, the optimal order (under some shape restrictions) for $\varepsilon_n$ is $n^{-1/2d}$. The estimator (4) (and the corresponding boundary estimator $T_n$) is obtained for several values of the smoothing parameter $\varepsilon_n$. The idea is to check the sensitivity of the estimation process with respect to changes in the value of $\varepsilon_n$. Alternatively, some procedure (cross-validation, bootstrap-based choice) for the optimal selection of the smoothing parameter could be used. However, in real applications (see Section 4.3 below) an optimal choice would be not so crucial when the procedure is used to establish comparisons between several images. In the case of a bi-color digitized image the calculation of $\mu(T_n)$ (and that of the area that appears in the denominator of the CI) is made by a simple count of the corresponding activated (black) pixels.

Note that, in practice, the first stage could be omitted as, strictly speaking, only the randomly selected points need to be classified. An interesting open problem in this regard would be to consider a more realistic model incorporating the classification error in the "red" or "green" areas, $G$ and $R$ (see Section 2).

4.3. *Results.* In the example of Figure 1 we have considered two sample sizes, $n = 50{,}000$ and $n = 100{,}000$. The results are summarized in Table 1.

The choices of the values $\varepsilon_n$ are of type $C_k n^{-1/4}$, where the constants $C_k$, for $k = 1, 2, 3$, are taken in order to consider small perturbations around the reference value $n^{-1/4}$. In the case $n = 100{,}000$ we have $n^{-1/4} = 0.0562$, so we decided to take $C_1$, $C_2$ and $C_3$ in order to get "exact" values (0.05, 0.02 and 0.01) for the smoothing parameter $C_k \varepsilon_n$. This entails that $C_1 = 0.8897$, $C_2 = 0.3559$ and $C_3 = 0.1779$ and we have kept these constants for the case $n = 50{,}000$.

The output in Table 1 indicates that the CI value is about 5.4. Clearly, the values (3.61, 3.96) obtained for the largest choices of $\varepsilon_n$ correspond to oversmoothed estimations; recall that the CI for a circle is 3.5449. This is apparent from the image of Figure 2, which shows the estimate $T_n$ of the infarct boundary for the case $n = 50{,}000$, $\varepsilon_n = 0.05$.

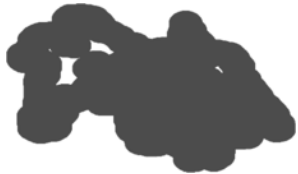

Fig. 2. *Oversmoothed boundary estimation of the infarct area in Figure* 1.

A remarkable fact in the results is their small variability. This means that, in practice, we can use a given (not necessarily optimal) choice of $\varepsilon_n$ to perform comparisons between different images. Maybe the true CI's are estimated with an appreciable bias, but this is, by far, the main source of error. Thus, the estimated CI's would allow us to get an assessment of the relative importance of the different cases from the point of view of infarct geometry, and the value of $\varepsilon_n$ corresponds, in some sense, to the resolution level employed in the procedure.

It is also worthwhile to observe that due to the presence of $\varepsilon_n$ in both the numerator and denominator of (4), the variability of this estimator is not a monotone function of $\varepsilon_n$. This is in contrast to the usual behavior of nonparametric estimators (e.g., kernel density estimators).

The estimation CI $\simeq 5.4$ suggests a rather negative diagnostic for the infarct shown in Figure 1. For example, in [11] the "infarct geometry" of a control group of eight infarcted pigs was studied and compared with that of another treatment group of eight individuals, also suffering a miochardial infarct but receiving a drug called 2,3-butanedione monoxime. The values found for the CI in the control and the treatment group are $7.7 \pm 0.2$ and $9.4 \pm 0.7$, respectively, which suggests a much better prognosis than that in our example.

In the case of the digital images, the choice of the smoothing parameter $\varepsilon_n$ is obviously limited by the pixel size. In our case, each side of the square $[0,1] \times [0,1]$ was divided into 600 square pixels so that the minimum effective choice of $\varepsilon_n$ would be $1/600 = 0.0017$.

TABLE 1
*Average values and standard deviations along* 100 *replications of the CI estimation for the infarct area in Figure* 1

| Sample size | $n = 50{,}000$ | | | $n = 100{,}000$ | | |
|---|---|---|---|---|---|---|
| $\varepsilon_n$ | **0.0119** | **0.0238** | **0.0595** | **0.01** | **0.02** | **0.05** |
| Mean | 5.2080 | 5.1265 | 3.6104 | 5.7257 | 5.53 | 3.96 |
| Standard deviation | 0.0042 | 0.00342 | 0.0129 | 0.0294 | 0.0213 | 0.0099 |

On the other hand, the large sample sizes ($n = 50{,}000, 100{,}000$) used in the study suggest the idea of using all the available pixels (360,000 in this example). The practical implications of such an "exhaustive method" are analyzed in some detail below [see paragraphs (g) and (h) in Section 5.2], in connection with the simulation example considered there, where the true value of the boundary length is exactly known.

The relative simplicity of the proposed method suggests the possibility of generalizations to multicolor higher-dimensional images; these could appear in the context of magnetic resonance explorations where very precise determinations are obtained for different magnitudes as the pH or the ATP (which measures the energy cell status).

**5. Simulations.** The estimator (4) is designed for cases where only incomplete information (given by "natural" sampling points on both sides of the border) is available. In this sense, the proposed method can be seen as a refined version of the nonparametric method for estimating boundaries discussed in [7]. The requirement of two samples (inside and outside the set) can be formalized with different models, but seems to be unavoidable in order to estimate the surface measure, unless we are willing to impose strong assumptions on the shape of $G$. On the other hand, the estimator (4) can be based on Monte Carlo (artificial) samples, to be used in contexts not directly related to image analysis, just as a stochastic device for approximating the length of a closed curve or the surface area of a body in $\mathbb{R}^3$.

As an example, we have considered the so-called Tschirnhausen Cubic (also known as Catalan's trisectrix and l'Hospital's cubic), a plane curve whose polar equations are

$$r = a\sec^3(\theta/3), \qquad \text{for } \theta \in (0, \pi),$$
$$r = a\sec^3((2\pi - \theta)/3), \qquad \text{for } \theta \in (\pi, 2\pi).$$

The reason for choosing this curve is the existence of closed simple expressions for both the length ($L_0 = 12a\sqrt{3}$) of the loop and the area inside ($A = 72a^2\sqrt{3}/5$). We have used our estimation method in order to approximate $L_0$ and $A$ in the case $a = 1$ (see Figure 3), so the target values are $L_0 = 20.7846$ and $A = 24.9415$. The random samples, with sizes $n = 30{,}000$ and $n = 10{,}000$, are drawn in the square $[-9, 2] \times [-5.5, 5.5]$, which fully includes the Tschirnhausen loop.

Before discussing the simulation experiment and output, we should consider a practical issue regarding the effective calculation of the estimator.

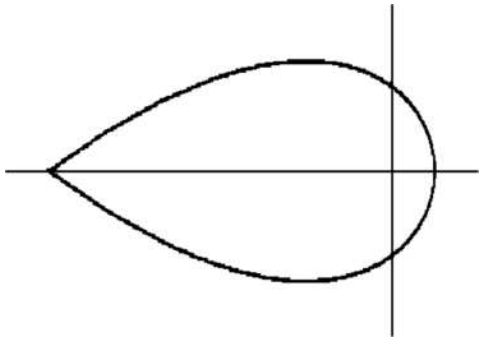

Fig. 3. *The Tschirnhausen Cubic.*

5.1. *Monte Carlo approximation of the estimator.* The estimator (4) (and the corresponding boundary estimator $T_n$) must be computed for every choice of $\varepsilon_n$ considered. An important practical problem is that the direct computation of $\mu(T_n)$ is not an easy task. However, it can be approximated easily, with an arbitrary precision level, by using again the Monte Carlo method. Let $Z_1^*, \ldots, Z_B^*$ be a random sample (independent of $Z_1, \ldots, Z_n$) from the uniform distribution on the unit square $[0,1]^d$. Since, with probability one, $T_n \subset [0,1]^d$ for $n$ large enough, we have that $\mu(T_n) = P(Z_1^* \in T_n)$ and therefore, for $B$ large,

$$\mu_B(T_n) = \frac{\sum_{i=1}^B \mathbb{I}_{\{Z_i^* \in T_n\}}}{B}$$

should be a good approximation of $\mu(T_n)$. Note that it is very easy to check when a point $Z^*$ belongs to $T_n$. This Monte Carlo method provides an approximate evaluation for $L_n$,

$$L_{n,B}^* = \frac{\mu_B(T_n)}{2\varepsilon_n}.$$

An interesting question in order to apply the proposed method is how to pick $B$ (as a function of $n$) to ensure that $L_{n,B}^*$ is a consistent estimator of $L_0$. The next theorem gives an answer to this question. The proof is given in Section 6.

THEOREM 3. *Besides the hypothesis of Theorem 1, let us assume that*

(12) $$\frac{B\varepsilon_n}{\log n} \to \infty.$$

*Then $L_{n,B}^* \to L_0$, a.s.*

5.2. *Simulation output.* The results of our simulation study are summarized in Table 2. The estimator $L_n$ has been evaluated for 500 samples of sizes $n = 30{,}000$ and $n = 10{,}000$. The resampling parameter $B$, used in the Monte Carlo approximations of $\mu(T_n)$, was $B = 1500$ in all cases considered. The output in Tables 2 (for $n = 30{,}000$) and 3 (for $n = 10{,}000$) provides the average, standard deviation and median of $L_n$ computed from the 500 replications, for different values of $\varepsilon_n$. The output is obtained using the same simulated samples for each value of $\varepsilon_n$. Thus the usual Monte Carlo area estimate, which does not depend on a smoothing parameter, is the same in all cases. The average, standard deviation and median obtained for this area estimator are respectively 24.9196, 0.4458 and 24.9125 for $n = 30{,}000$ and 24.9485, 0.7842 and 24.9889 for $n = 10{,}000$.

Some direct conclusions can be drawn from these results:

(a) The true value $L_0 = 20.7846$ is systematically underestimated with a relative error about 4.7% (in the case $n = 30{,}000$) and 8.1% (for $n = 10{,}000$). The gain obtained by increasing the sample size is mostly apparent in the bias. The average (over all values of $\varepsilon_n$) of the average output is 19.0890 for $n = 10{,}000$ and 19.7900 for $n = 30{,}000$.

TABLE 2
*Average, standard deviation and median of $L_n$ computed over 500 replications with $n = 30{,}000$*

| $\varepsilon_n$ | 0.76 | 0.78 | 0.80 | 0.82 | 0.84 | 0.86 | 0.88 |
|---|---|---|---|---|---|---|---|
| Average | 19.7301 | 19.7416 | 19.7621 | 19.7644 | 19.7918 | 19.7918 | 19.7859 |
| Std. deviation | 1.3940 | 1.3935 | 1.3793 | 1.3448 | 1.3470 | 1.3200 | 1.3072 |
| Median | 19.7548 | 19.7920 | 19.8274 | 19.7576 | 19.8930 | 19.8249 | 19.8081 |
| $\varepsilon_n$ | 0.9 | 0.92 | 0.94 | 0.96 | 0.98 | 1.0 | 1.2 |
| Average | 19.7901 | 19.7949 | 19.8109 | 19.8208 | 19.8290 | 19.8230 | 19.8237 |
| Std. deviation | 1.2952 | 1.2917 | 1.2636 | 1.2331 | 1.2159 | 1.2031 | 1.0666 |
| Median | 19.8863 | 19.9150 | 19.8522 | 19.8804 | 19.8209 | 19.8804 | 19.8627 |

TABLE 3
*Average, standard deviation and median of $L_n$ computed over 500 replications with $n = 10{,}000$*

| $\varepsilon_n$ | 0.76 | 0.78 | 0.80 | 0.82 | 0.84 | 0.86 | 0.88 |
|---|---|---|---|---|---|---|---|
| Average | 19.0026 | 18.8594 | 18.9512 | 18.9370 | 18.9507 | 19.0806 | 19.1083 |
| Std. deviation | 1.3908 | 1.3586 | 1.3260 | 1.2627 | 1.3075 | 1.3398 | 1.2467 |
| Median | 18.9736 | 18.9221 | 18.9791 | 18.9300 | 18.9842 | 19.0359 | 19.0852 |
| $\varepsilon_n$ | 0.9 | 0.92 | 0.94 | 0.96 | 0.98 | 1.0 | 1.2 |
| Average | 19.0492 | 19.1409 | 19.1408 | 19.2057 | 19.2134 | 19.2384 | 19.3679 |
| Std. deviation | 1.2956 | 1.1936 | 1.1599 | 1.2460 | 1.2157 | 1.1777 | 1.0394 |
| Median | 19.0381 | 19.1774 | 19.1303 | 19.1735 | 19.2150 | 19.2548 | 19.3679 |

(b) The simulation output shows a considerable stability with respect to the values of the smoothing parameter $\varepsilon_n$. This stability remains even for other smaller values of $\varepsilon_n$ (not included in the tables) that we have checked. For example, whereas the average of the average values of $L_n$ obtained from the 14 choices of $\varepsilon_n$ included in Table 2 is 19.79, the corresponding average for the other five equispaced values of $\varepsilon_n$, ranging from 0.64 to 0.72, is 19.5985.

(c) The sampling distributions are almost symmetric, with the median very close to the mean in all cases.

(d) There is a slight, but consistent, decline of the variability around the mean as the smoothing parameter increases.

(e) As could be predicted, the variability results tend to improve, at the expense of some additional computational burden, by increasing the value of the resampling parameter $B$. For example, the output for the average, standard deviation and median of $L_n$ with $\varepsilon_n = 0.92$, $n = 30{,}000$ and $B = 2000$ is 19.8341, 1.0790 and 19.8458, respectively. For $n = 10{,}000$, with the same value of $\varepsilon_n$, the corresponding output for $B = 3000$ is 19.2229, 0.8490 and 19.1774. These results account for the small changes in the variability of $L_n$ from $n = 10{,}000$ to $n = 30{,}000$. They suggest that, for these sample size magnitudes, most variability is due to the Monte Carlo approximation stage of the numerator $\mu(T_n)$ in (4), controlled by the parameter $B$. The value $B = 1500$, used in the simulations of Tables 2 and 3, should be considered as a first computationally affordable choice, suitable for this preliminary study.

(f) The plots of the density estimators obtained from the values of $L_n$ suggest that the sampling distribution is, for all the considered choices of $\varepsilon_n$, very close to normality. As a consequence, an interesting open problem would be to establish the asymptotic normality of $L_n$. However, the proof seems far from trivial in view of the special structure of $L_n$.

(g) In this example we have implemented our method in a case where an exact equation for the border is known. So no digitization process is involved. In practice, most real black-and-white images come in a digitized version. In mathematical terms this amounts to replacing the original image $G$ by a finite union $G^h$ of square pixels with sides of fixed length $h$, parallel to the coordinate axes. In such situations one could think of exactly measuring the border length $L^h$ of the "digital boundary" $\partial G^h$. This is just the number of pixel sides separating regions of different colors. This is computationally feasible and avoids the use of any smoothing parameter. However, it is not difficult to see that this direct exhaustive procedure will fail, as $L^h$ cannot converge to $L$ when $h$ tends to 0. For example, if $\partial G$ includes a segment $A$ inside the diagonal $x = y$, the length of $A$ will be overestimated by a factor $\sqrt{2}$ when $G$ is approximated by $G^h$. This is empirically confirmed in our case:

If we replace the region $G$ inside the Tschirnhausen Cubic by a digitized version, obtained by dividing the "frame square" $[-9, 2] \times [-5.5, 5.5]$ into $300 \times 300$ pixels, the direct exhaustive method gives an estimation $L^h = 25.97$ for the true value $L_0 = 20.78$. Our method, with $n = 10{,}000$, provides much more acceptable estimation around 19.7 (see Table 3). The use of a more precise digitization does not improve things (in fact, it reveals a lack of consistency in the exhaustive procedure). For example, a $600 \times 600$ digitization leads to $L^h \simeq 26.5$, and with $1024 \times 1024$ pixels we get $L^h \simeq 28.1875$.

The exhaustive method could also been implemented in an indirect version, based on measuring areas. The boundary length could be estimated by $\text{Area}(G^{hb})/2h$, where $G^{hb}$ denotes the union of all "boundary pixels" in $G^h$. This also fails: estimation for the $300 \times 300$ and $600 \times 600$ digitizations gives respectively 19.36 and 19.32. Note that, in fact, this procedure uses implicitly a smoothing parameter (the pixel side length $h$). The failure should be interpreted as a phenomenon of undersmoothing; see the comment about the bias after (10) and (11).

(h) The use of all the available pixels is still a possibility, although, in view of the previous comment it should always be done with an appropriate amount of smoothing, along the lines indicated above. Although this exhaustive procedure "with smoothing" is feasible in many cases, it is not advisable in general, due to its lack of robustness against the "noise" (in the form of disperse error pixels not belonging to the image). By contrast, the method based on random samples will automatically ignore (with high probability) the possible disperse noise, at the expense of higher variability. We have checked this by randomly adding four patches of noise, in the form of circular clusters (with radii 0.25) of black pixels, within the square $[-9, 2] \times [-5.5, 5.5]$, where the loop of the Tschirnhausen Cubic is included. In the worst case (when the four noise patches fall on the white background, outside the black image), the amount of noise added to the image represents less than 1% of the total number of pixels. The presence of the noise turned out to have a devastating effect in the exhaustive method with smoothing: The average length obtained with this method for 500 of such noisy images is 24.92 (standard deviation 1.63), whereas the random method applied with a sample size $n = 5000$ and $\varepsilon_n = 0.94$ gave an average of 21.07 (standard deviation 0.9992). Curiously enough, the results for the latter method (recall that the true value is 20.78) are even better than those obtained in the case with no noise since the noise tends to increase the boundary length, thus partially correcting the inherent underestimation bias.

## 6. Proofs.

PROOF OF THEOREM 1. The result is a consequence of the following two claims.

STATEMENT 1. *With probability one, $T_n \subset B(T, \varepsilon_n)$.*

STATEMENT 2. *For any $0 < \alpha < 1$, we have eventually, with probability one,*

$$B(T, \varepsilon'_n) \subset T_n,$$

*where $\varepsilon'_n = \alpha \varepsilon_n, 0 < \alpha < 1$.*

PROOF OF STATEMENT 1. For any $z \in T_n$, we have that (with probability one) $B(z, \varepsilon_n)$ meets $G$ and its complementary $R$. Therefore, $B(z, \varepsilon_n)$ meets the boundary of $G$, $T$, which means that $z$ belongs to $B(T, \varepsilon_n)$. This concludes the proof of Statement 1. □

PROOF OF STATEMENT 2. By the Borel–Cantelli lemma, it is sufficient to show that

$$\sum_{n=1}^{\infty} P(B(T, \varepsilon'_n) \nsubseteq T_n) < \infty.$$

However,

(13) $$P(B(T, \varepsilon'_n) \nsubseteq T_n) \leq P(\exists z \in B(T, \varepsilon'_n) : G_z(\varepsilon_n) = 0)$$
$$+ P(\exists z \in B(T, \varepsilon'_n) : R_z(\varepsilon_n) = 0).$$

Now, we try to find an upper bound for the first probability on the right-hand side. The other probability can be bounded by a similar argument.

For any $z \in B(T, \varepsilon'_n)$, there is an $t \in T$ for which $B(t, \beta_n) \subset B(z, \varepsilon_n)$, where $\beta_n = (1 - \alpha)\varepsilon_n$. Therefore,

$$P(\exists z \in B(T, \varepsilon'_n) : G_z(\varepsilon_n) = 0) \leq P(\exists t \in T : G_t(\beta_n) = 0).$$

Let $T(\beta_n)$ be a set [with cardinality $D(\beta_n)$] of ball centres corresponding to a minimal covering of $T$ by balls of radius $\beta_n/2$. So we consider a class $\{B(s, \beta_n/2) : s \in T(\beta_n) \subset T\}$ such that

$$T \subset \bigcup_{s \in T(\beta_n)} B\left(s, \frac{\beta_n}{2}\right).$$

Since $\{\exists t \in T : G_t(\beta_n) = 0\} \subset \{\exists s \in T(\beta_n) : G_s(\beta_n/2) = 0\}$, we have

$$P(\exists t \in T : G_t(\beta_n) = 0) \leq P\left(\exists s \in T(\beta_n) : G_s\left(\frac{\beta_n}{2}\right) = 0\right)$$
$$\leq \sum_{s \in T(\beta_n)} P\left(G_s\left(\frac{\beta_n}{2}\right) = 0\right)$$

$$= \sum_{s \in T(\beta_n)} \left(1 - p_X\left(s, \frac{\beta_n}{2}\right)\right)^n$$

$$\leq \sum_{s \in T(\beta_n)} \exp\left\{-n p_X\left(s, \frac{\beta_n}{2}\right)\right\},$$

where in the last inequality we have used the fact that $1 - x \leq e^{-x}$ for $0 \leq x \leq 1$. The right-hand side of the above inequality can easily be bounded since, from the standardness hypothesis, for $n$ large enough,

$$p_X\left(s, \frac{\beta_n}{2}\right) \geq C \omega_d \mu(G) \frac{\beta_n^d}{2^d} = K_1 \varepsilon_n^d,$$

where $\omega_d = \mu(B(0,1))$ and $K_1$ is a constant which depends on the dimension $d$, $\alpha$, $\mu(G)$ and $C$. Therefore,

$$P(\exists z \in B(T, \varepsilon_n') : G_z(\varepsilon_n) = 0) \leq D(\beta_n) \exp\{-K_1 \varepsilon_n^d\}.$$

Now, in order to bound the function $D(\varepsilon)$, recall that it represents the cardinality of a minimal covering $\mathcal{C}(\varepsilon/2)$ of $T$ by balls of radii $\varepsilon/2$. This entails (e.g. [14], page 78) that there exists a family of $D(\varepsilon)$ disjoint balls with radii $\varepsilon/4$ and centres at points of $T$. Then the sum of their measures must be smaller than $\mu(B(T, \varepsilon/4))$. Hence,

$$D(\varepsilon)(\varepsilon/4)^d \omega_d \leq \mu(B(T, \varepsilon/4)).$$

Since $L(\varepsilon) \to L_0$, we get for $\varepsilon$ small enough, $D(\varepsilon) \leq A\varepsilon^{1-d}$ for some constant $A$. Therefore,

$$P(\exists z \in B(T, \varepsilon_n') : G_z(\varepsilon_n) = 0) \leq K_2 \varepsilon_n^{1-d} \exp(-K_1 n \varepsilon_n^d),$$

where $K_2 = (1-\alpha)^{1-d} A$. The condition $n\varepsilon_n^d/\log n \to \infty$ ensures the convergence of the series $\sum_{n=1}^{\infty} \varepsilon_n^{1-d} \exp(-K_1 n \varepsilon_n^d)$. The other probability in (13) can be bounded in a similar way. Note that the obvious inequality $D(\varepsilon) \leq A\varepsilon^{-d}$ would also suffice for the purpose of convergence, but the above simple argument provides a sharper bound for the probabilities. This concludes the proof of Statement 2. □

Now the proof of Theorem 1 is a straightforward consequence of Statements 1 and 2. Indeed, we have that, with probability one,

$$\alpha L_0 = \lim_n \frac{\mu(B(T, \varepsilon_n'))}{2\varepsilon_n} \leq \liminf_n L_n \leq \limsup_n L_n \leq \lim_n \frac{\mu(B(T, \varepsilon_n))}{2\varepsilon_n} = L_0.$$

This holds for any $\alpha \in (0,1)$ and therefore, the conclusion of the theorem follows. □

PROOF OF THEOREM 2. The expected value of $L_n$ can be written as

$$E(L_n) = \frac{E(\mu(T_n))}{2\varepsilon_n} = \frac{1}{2\varepsilon_n} E\left(\int \mathbb{I}_{\{z \in T_n\}} \mu(dz)\right) = \frac{1}{2\varepsilon_n} \int E(\mathbb{I}_{\{z \in T_n\}}) \mu(dz)$$

$$= \frac{1}{2\varepsilon_n} \int P(z \in T_n) \mu(dz) = \frac{1}{2\varepsilon_n} \int_{B(T,\varepsilon_n)} P(z \in T_n) \mu(dz),$$

since, with probability one, $T_n \subset B(T, \varepsilon_n)$. It is clear that

(14) $$P(z \notin T_n) \leq P(G_z(\varepsilon_n) = 0) + P(R_z(\varepsilon_n) = 0).$$

Remember that $G_z(\varepsilon_n)$ has a binomial distribution with parameters $n$ and $p_X(z, \varepsilon_n)$. Therefore,

$$P(G_z(\varepsilon_n) = 0) = (1 - p_X(z, \varepsilon_n))^n \leq \exp\{-np_X(z, \varepsilon_n)\}.$$

Let $P_T z \in T$ be the projection of $z$ onto $T$. Since, for any $z \in B(T, \varepsilon_n)$,

$$B(P_T z, \varepsilon_n - d(z, T)) \subset B(z, \varepsilon_n),$$

using condition (a) of Theorem 1, we have that, for $\varepsilon_n$ small enough,

$$P_X(B(z, \varepsilon_n)) \geq C\omega_d(\varepsilon_n - d(z, T))^d.$$

Hence,

$$P(G_z = 0) \leq \exp\{-K_1 n(\varepsilon_n - d(z, T))^d\},$$

where $K_1$ is a positive constant which depends only on $\mu(G)$, $C$ and the dimension $d$. Similarly, we have that $P(R_z = 0) \leq \exp\{-K_2 n(\varepsilon_n - d(z, T))^d\}$, for a positive constant $K_2$ which depends only on $\mu(R)$, $C$ and $d$. Using these bounds and (14), we get

$$P(z \in T_n) \geq 1 - 2\exp\{-Kn(\varepsilon_n - d(z, T))^d\},$$

where $K = \min(K_1, K_2)$. Thus, we have that

$$E(L_n) = \frac{1}{2\varepsilon_n} \int_{B(T,\varepsilon_n)} P(z \in T_n) \, dz$$

$$\geq \frac{1}{2\varepsilon_n} \int_{B(T,\varepsilon_n)} (1 - 2\exp\{-Kn(\varepsilon_n - d(z, T))^d\}) \, dz$$

$$= L(\varepsilon_n) - \frac{1}{\varepsilon_n} \int_{B(T,\varepsilon_n)} \exp\{-Kn(\varepsilon_n - d(z, T))^d\} \, dz = L(\varepsilon_n) - I_n,$$

with

$$I_n = \frac{1}{\varepsilon_n} \int_{B(T,\varepsilon_n)} g_n(d(z, T)) \, dz,$$

where $g_n(w) = \exp\{-Kn(\varepsilon_n - w)^d\}$. By the change of variable formula, we have that

$$(15) \qquad I_n = \frac{1}{\varepsilon_n} \int_0^{\varepsilon_n} g_n(w) F(dw),$$

where $F(w) = \mu(\{z : d(z,T) \le w\}) = \mu(B(T,w))$. By the assumption made on the continuous differentiability of $F$ at 0 and the existence and finiteness of the Minkowski content, we have $F'(0) = 2L_0$ so that, for $w$ small enough, $F'(\omega) \le 3L_0$. Finally, for $n$ large enough,

$$I_n \le \frac{3L_0}{\varepsilon_n} \int_0^{\varepsilon_n} \exp\{-Knt^d\}\, dt = \frac{3L_0}{\varepsilon_n} \int_0^{Kn\varepsilon_n^d} \exp(-u) \frac{1}{d(Kn)^{1/d}} u^{-(d-1)/d}\, du$$

$$\le \frac{3L_0}{dK^{1/d}(\varepsilon_n^d n)^{1/d}} \int_0^\infty \exp(-u) u^{-(d-1)/d}\, du = \frac{A}{(\varepsilon_n^d n)^{1/d}},$$

where in the first inequality we have applied in (15) the change of variable $t = \varepsilon_n - w$ and then (for the first equality) $u = Knt^d$. □

PROOF OF COROLLARY 1. The bound (9) for the $L_1$-error follows as a direct consequence of the bounds (6)–(8) together with the triangle inequality. Now, the conclusion (a) follows from (8) and the definition of $L_0$.

To show (b), note that the optimal convergence order for the bound (9) is obtained by making equal the convergence orders of both terms on the right-hand side. Under the smoothness conditions mentioned in Section 3.2, we have $|L(\varepsilon_n) - L_0| = O(\varepsilon_n)$ (see [10], Theorem 5.6). Thus, from (8), the optimal order for the bound (9) is $O(n^{-1/2d})$, which is attained for $\varepsilon_n = n^{-1/2d}$. □

PROOF OF THEOREM 3. Clearly, it is enough to show that $L_{n,B}^* - L_n \to 0$, a.s. This can be proved showing that, for any $\rho > 0$,

$$(16) \qquad \sum_n P(|L_{n,B}^* - L_n| > \rho) < \infty.$$

This is not hard to do because, given $Z_1, \ldots, Z_n$, $L_{n,B}^*$ has (essentially) a binomial distribution with mean $L_n$ and, therefore, we can use a concentration inequality to control the size of its tail. Indeed,

$$P(|L_{n,B}^* - L_n| > \rho)$$
$$= E(P(|L_{n,B}^* - L_n| > \rho | Z_1, \ldots, Z_n))$$
$$= E(P(|\mu_B(T_n) - \mu(T_n)| > 2\rho\varepsilon_n | Z_1, \ldots, Z_n))$$
$$\le E\left(2\exp\left\{-\frac{4\rho^2 \varepsilon_n^2 B}{2\mu(T_n)(1-\mu(T_n)) + (4/3)\rho\varepsilon_n}\right\}\right),$$

where in the last step we have used Bernstein's inequality. It is not hard to bound this last quantity because $\mu(T_n)$ goes to zero (with probability one) as fast as $\varepsilon_n$ when $n$ tends to infinity. To see this, note that in Theorem 1 we proved that (with probability one) $T_n \subset B(T, \varepsilon_n)$ and, therefore, $\mu(T_n) \leq \mu(B(T, \varepsilon_n))$. Since $L(\varepsilon_n) \to L_0$, we have that, for $n$ large enough, $\mu(B(T, \varepsilon_n)) \leq 4L_0\varepsilon_n$. So, for $n$ large enough,

$$E\left(2\exp\left\{-\frac{4\rho^2\varepsilon_n^2 B}{2\mu(T_n)(1-\mu(T_n))+(4/3)\rho\varepsilon_n}\right\}\right)$$
$$\leq E\left(2\exp\left\{-\frac{4\rho^2\varepsilon_n^2 B}{8L_0\varepsilon_n+(4/3)\rho\varepsilon_n}\right\}\right)$$
$$= 2\exp\{-K_{\rho,L_0}\varepsilon_n B\},$$

where $K_{\rho,L_0}$ is a (positive) constant. Obviously, (12) ensures that, for any $\rho > 0$,

$$\sum_n \exp\{-K_{\rho,L_0}\varepsilon_n B\} < \infty,$$

and, therefore, (16) holds. This concludes the proof of the theorem. $\square$

**Acknowledgments.** We are most grateful to Dr. David García-Dorado (Servei de Cardiologia, Hospital Vall d'Hebron, Barcelona) who provided us with the material for the cardiology case-study and explained to us the required medical concepts. We are also indebted to Professor Jesús Gonzalo (Departamento de Matemáticas, Universidad Autónoma de Madrid) for his very helpful geometrical insights on the Minkowski content notion. The third author also wishes to thank the University of Vigo (Spain), where he carried out part of his work on this paper. The constructive comments of two referees are gratefully acknowledged.

A. CUEVAS  
DEPARTAMENTO DE MATEMÁTICAS  
UNIVERSIDAD AUTÓNOMA DE MADRID  
28049 MADRID  
SPAIN  
E-MAIL: antonio.cuevas@uam.es

R. FRAIMAN  
DEPARTAMENTO DE MATEMÁTICA  
UNIVERSIDAD DE SAN ANDRÉS  
ARGENTINA  
E-MAIL: rfraiman@udesa.edu.ar

A. RODRÍGUEZ-CASAL  
DEPARTAMENTO DE ESTATÍSTICA E INV. OPERATIVA  
UNIVERSIDAD DE SANTIAGO DE COMPOSTELA  
SPAIN  
E-MAIL: alrodcas@usc.es